\documentclass[english,draft]{article}
\usepackage[english]{babel}
\usepackage{amsmath}
\usepackage{amsthm}
\usepackage{amssymb}
\usepackage{color}
\usepackage{textcomp}

\newtheorem{theorem}{Theorem}[]
\newtheorem{lemma}{Lemma}[]

\newtheorem{problem}{\textbf{Problem}}[]

\newtheorem{statement}{Statement}[]
\newtheorem{corollary}{Corollary}[]
\newtheorem{definition}{Definition}[]

\mathversion{bold}
\begin{document}
\title{\textbf{The concrete theory of numbers :\\ New Mersenne conjectures. Simplicity and other wonderful properties of numbers $L(n) = 2^{2n}\pm2^n\pm1$.}}
\author{\textbf{Boris\, V. Tarasov}\thanks{Tarasov Boris V.
The concrete theory of numbers: New Mersenne conjectures. Simplicity and other wonderful properties of numbers $L(n) = 2^{2n}\pm2^n\pm1$.
MSC 11A51, MSC 11B83.
\textbf{\textcopyright 2008} \textbf{Tarasov Boris V.}
}\\
}
\maketitle
\begin{abstract}
New Mersenne conjectures. The problems of simplicity, common prime divisors and
free from squares of numbers $L(n) = 2^{2n}\pm2^n\pm1$
are investigated. Wonderful formulas $gcd $ for numbers $L (n) $ and numbers repunit are proved.
\end{abstract}

{\centering\section{Introduction}\par}

In present work we consider sequences of integers of the following kind :
\begin{gather}
\textcolor[rgb]{0.00,0.00,0.50}{L_{1}(n) = 2^{2n} + 2^{n} + 1},\tag{$L_{1}$}\label{e:L1}\\
\textcolor[rgb]{0.00,0.00,0.50}{L_{2}(n) = 2^{2n} + 2^{n} - 1},\tag{$L_{2}$}\label{e:L2}\\
\textcolor[rgb]{0.00,0.00,0.50}{L_{3}(n) = 2^{2n} - 2^{n} + 1},\tag{$L_{3}$}\label{e:L3}\\
\textcolor[rgb]{0.00,0.00,0.50}{L_{4}(n) = 2^{2n} - 2^{n} - 1},\tag{$L_{4}$}\label{e:L4}
\end{gather}where $n\geq1$ is integer.\par
For the numerical sequence, being the union  of numerical sequences $L_{1}(n)$, $L_{2}(n)$,
$L_{3}(n)$, $L_{4}(n)$, we use a designation
\begin{equation}\textcolor[rgb]{0.00,0.00,0.50}{L(n) = 2^{2n} \pm 2^{n} \pm 1},\label{e:L}\end{equation}
where $n\geq1$ is integer.\par
The author is interested to research the new Mersenne conjectures, concerning numbers $L(n)$. The reviews concerning Mersenne numbers
 and new Mersenne conjectures are available here
\cite{Graham, WeissteinF, WeissteinM, Mersenne, Bateman}.\par
\par Will we use the following designations further :
\par $(a,b)=gcd(a,b)$~is the greatest common divider of integers $a>0$, $b>0$.
\par $p,q$ are odd prime numbers.
\par $m \bot n \Longleftrightarrow m, n -$ are integers and $gcd(m,n)=1$(see\cite{Graham}).
\par If it is not stipulated specially, the integer positive numbers are considered.\par

We are interested by the following questions concerning numbers $L(n)$ : question of simplicity of numbers, question
of the common divisors and question of freedom from squares.\par
For numbers $L(n)$ two general simple statements are fair. These statements represent trivial consequences of the small Fermat's theorem ,
submitted by the following comparison,\\
\begin{center} \textcolor[rgb]{0.00,0.00,0.50}{$n^{p^{k}}\equiv n^{p^{k-1}}(mod\,p^{k})$} for $k>0$(see\cite{Graham}) \end{center}
or chain of comparisons \\
\begin{center} \textcolor[rgb]{0.00,0.00,0.50}{$n^{p^{N+k}}\equiv n^{p^{N-1+k}}\equiv\ldots\equiv n^{p^{k}}\equiv n^{p^{k-1}}(mod\,p^{k})$}
for $k>0$, $N\geq0$,\end{center} $n$ of integers.
\begin{statement}[]\label{S:Statement1}
If $L(l)\equiv0(mod\,p)$, where $l>0$ is integer, $p$ is prime number, then
\begin{displaymath}
\textcolor[rgb]{0.00,0.00,0.50}{L((p-1)\cdot k+l)\equiv0(mod\,p)},
\end{displaymath}
where $k\geq0$ is any integer.
\end{statement}
\begin{statement}[]\label{S:Statemen2}
If $L(l)\equiv0(mod\,p^{t})$, where $t>0$, $l>0$ is integer, $p$ is prime number, then
\begin{displaymath}
\textcolor[rgb]{0.00,0.00,0.50}{L(p^{N+t}-p^{t-1}+l)\equiv0(mod\,p^{t})},
\end{displaymath}
where $N\geq0$ is any integer.
\end{statement}

{\centering\section{Numbers $L_{1}(n) = 2^{2n} + 2^{n} + 1$}\par}

\paragraph{\textbf{\textsection 1}.} Let's bring the simple statements concerning numbers $L_{1}(n)$.
\begin{lemma}[]\label{L:L1} For numbers $L_{1}(n) = 2^{2n} + 2^{n} + 1$ following statements are\\ fair :\par
(1) The first prime numbers $L_{1}(n)$ correspond to $n=1, 3, 9$.\\ $L_{1}(1)\,=\,7$,
$L_{1}(3)\,=\,73$,
$L_{1}(9)\,=\,262657$.\par
(2) If $n\geq 2$ is an even number, then number
\begin{equation}\label{e:N1} \textcolor[rgb]{0.00,0.00,0.50}{L_{1}(n) \equiv 0(mod\,3)} \end{equation}
is composite.
If $n\geq 1$ is an odd number, then
\begin{equation}\label{e:N2} \textcolor[rgb]{0.00,0.00,0.50}{L_{1}(n) \equiv 1\not\equiv 0(mod\,3)} \end{equation}\par
(3) If $n>1$, $n\not\equiv 0(mod\,3)$, then number
\begin{equation}\label{e:N3} \textcolor[rgb]{0.00,0.00,0.50}{L_{1}(n) \equiv 0(mod\,7)} \end{equation}
is composite.\par
\end{lemma}
\begin{proof}
Validity of congruences \eqref{e:N1}-\eqref{e:N3} obviously \\follows from trivial comparisons
$2 \equiv -1(mod\,3)$, $2^3 \equiv 1(mod\,7)$.
\end{proof}
Thus, prime numbers $L_{1}(n)$ are probable only for $n = 3k$, where $k$ is odd number.
\begin{lemma}[]\label{L:L2} Let $k>1$ is an integer. Then there will be a prime number $q>3$, such that
number 2 on the module $q$ belongs to a index $3^k$.
\end{lemma}
\begin{proof}
Let's consider expression\\ $A = 2^{3^k}-1 = (2^{3^{k-1}})^3-1 = (2^{3^{k-1}}-1)\cdot((2^{3^{k-1}})^2+2^{3^{k-1}}+1)$.
\\Let $q$ is a prime number such that\\ $B = (2^{3^{k-1}})^2+2^{3^{k-1}}+1\equiv0(mod\,q)$. If $q=3$, then
$B \equiv 1\not\equiv 0(mod\,3)$. Hence $q>3$, $2^{3^k}-1\equiv0(mod\,q)$. Let's assume, that\\
$2^{3^l}-1\equiv0(mod\,q)$, where $k>l\geq0$. Let $d=k-1-l\geq0$, $3^d\geq1$.\\ Then $(2^{3^l})^{3^d}\equiv1(mod\,q)$,
$2^{3^{l+d}}\equiv1(mod\,q)$, $2^{3^{k-1}}\equiv1(mod\,q)$,\\
$B = (2^{3^{k-1}})^2+2^{3^{k-1}}+1\equiv3\equiv0(mod\,q)$, $q=3$. Have received the contradiction. It is obvious, that $d=3^k$
 is the least positive number, for which the comparison $2^{d}-1\equiv0(mod\,q)$ is feasible.
\end{proof}
\paragraph{\textbf{\textsection 2}.} In connection with the lemma 1 the research of prime divisors of numbers $L_{1}(n)$ is interesting, where $n=3^{k}t$, where $k\geq1$
is integer, $t\geq1$ is an odd number. Since only at such $n$ the numbers $L_{1}(n)$  can be \emph{suspicious} on prime numbers.
In the following theorem the interesting property $gcd$ for numbers $L_{1}(n)$ is proved.
\begin{theorem}[$gcd$ of numbers $L_{1}$]\label{T:T1}
Let $k\geq0$, $k_{1}\geq0$, $k_{2}\geq0$ be integers; $t_{1}\geq1$, $t_{2}\geq1$ are odd numbers, $t_{1}\not\equiv0(mod\,3)$, $t_{2}\not\equiv0(mod\,3)$.
Then the following statements are fair :
\begin{gather}
\label{T1:1}\textcolor[rgb]{0.00,0.00,0.50}{gcd(L_{1}(3^{k}t_{1}),L_{1}(3^{k}t_{2}))=L_{1}(3^{k}gcd(t_{1},t_{2}))}.\\
\label{T1:2}\textcolor[rgb]{0.00,0.00,0.50}{gcd(L_{1}(3^{k_{1}}t_{1}),L_{1}(3^{k_{2}}t_{2})) = 1}
\end{gather}
at $k_{1}\neq k_{2}$.
\end{theorem}
\begin{proof}
Let $t_{3}=gcd(t_{1},t_{2})$, $t_{1}=t_{3}d_{1}$, $t_{2}=t_{3}d_{2}$, where \\$(d_{1},d_{2})=1$, $d_{1}\geq1$, $d_{2}\geq1$ are odd numbers.\par
1) Let's prove equality \eqref{T1:1}. Let's consider the following formulas
\begin{eqnarray}
\textcolor[rgb]{0.00,0.00,0.50}{A=L_{1}(3^{k}t_{1})=2^{2\cdot 3^{k}t_{1}} + 2^{3^{k}t_{1}} + 1,}
&\textcolor[rgb]{0.00,0.00,0.50}{A(2^{3^{k}t_{1}}-1)=2^{3^{k+1}t_{3}d_{1}}-1},\notag\\
\textcolor[rgb]{0.00,0.00,0.50}{B=L_{1}(3^{k}t_{2})=2^{2\cdot 3^{k}t_{2}} + 2^{3^{k}t_{2}} + 1,}
&\textcolor[rgb]{0.00,0.00,0.50}{B(2^{3^{k}t_{2}}-1)=2^{3^{k+1}t_{3}d_{2}}-1},\label{T1:3}\\
\textcolor[rgb]{0.00,0.00,0.50}{C=L_{1}(3^{k}t_{3})=2^{2\cdot 3^{k}t_{3}} + 2^{3^{k}t_{3}} + 1,}
&\textcolor[rgb]{0.00,0.00,0.50}{C(2^{3^{k}t_{3}}-1)=2^{3^{k+1}t_{3}}-1}.\notag
\end{eqnarray}
Then the following formulas are fair
\begin{multline}
\textcolor[rgb]{0.00,0.00,0.50}{A(2^{3^{k}t_{1}}-1)=C(2^{3^{k}t_{3}}-1)[2^{3^{k+1}t_{3}(d_{1}-1)} + 2^{3^{k+1}t_{3}(d_{1}-2)} + \ldots}\\
\textcolor[rgb]{0.00,0.00,0.50}{\ldots + 2^{3^{k+1}t_{3}} + 1]},\label{T1:4}
\end{multline}
\begin{multline}
\textcolor[rgb]{0.00,0.00,0.50}{B(2^{3^{k}t_{2}}-1)=C(2^{3^{k}t_{3}}-1)[2^{3^{k+1}t_{3}(d_{2}-1)} + 2^{3^{k+1}t_{3}(d_{2}-2)} + \ldots}\\
\textcolor[rgb]{0.00,0.00,0.50}{\ldots + 2^{3^{k+1}t_{3}} + 1]}.\label{T1:5}
\end{multline}
Let $q>1$ is prime number such, that
\begin{equation}
\textcolor[rgb]{0.00,0.00,0.50}{2^{3^{k}t_{1}}-1\equiv 0(mod\,q),\ C\equiv 0(mod\,q)}.\label{T1:6}
\end{equation}
Then it follows from \eqref{T1:3} and \eqref{T1:6}, that
\[
\textcolor[rgb]{0.00,0.00,0.50}{2^{3^{k}t_{3}d_{1}}-1\equiv 0(mod\,q),\ 2^{3^{k+1}t_{3}}-1\equiv 0(mod\,q)}.
\]
Let's consider number $b=2^{3^{k}t_{3}}$. If $b\equiv 1(mod\,q)$, then it follows from \eqref{T1:3} and \eqref{T1:6}, that $C\equiv 3\equiv 0(mod\,q)$, i.e. $q=3$.
Since $t_{3}$ is odd number, then it follows from lemma 1, that $C\not\equiv 0(mod\,3)$. We have come to the contradiction.\par Thus, $b\not\equiv 1(mod\,q)$.
Then the number $b$ on the module $q$ belongs to index $l_{0}>1$, $b^{l_{0}}\equiv 1(mod\,q)$. As $b^{d_{1}}\equiv b^{3}\equiv 1(mod\,q)$, that
$d_{1}\equiv 3\equiv 0(mod\,l_{0})$, $d_{1}\equiv 0(mod\,3)$. Have received the contradiction.\par
We have proved, that $gcd(C,2^{3^{k}t_{1}}-1)=gcd(C,2^{3^{k}t_{2}}-1)=1$, hence, $gcd(A,B)\equiv 0(mod\,C)$. It is necessary to prove the opposite: if\\
$d\mid gcd(A,B)$, then $d\mid C$, where $d>1$ an integer.\par Let's assume, that there is an integer $d>1$ such, that
\begin{equation}
\textcolor[rgb]{0.00,0.00,0.50}{A\equiv B\equiv 0(mod\,d),\ but\ gcd(C,d)=1}.\label{T1:7}
\end{equation}
Then it follows from \eqref{T1:3} and \eqref{T1:7}, that
\[
\textcolor[rgb]{0.00,0.00,0.50}{2^{3^{k+1}t_{3}d_{1}}-1\equiv 0(mod\,d),\ 2^{3^{k+1}t_{3}d_{2}}-1\equiv 0(mod\,d)}.
\]
Let $l_{0}>1$ is an index, to which the number $2$ belongs on the module $d$, \\$2^{l_{0}}-1\equiv 0(mod\,d)$. Then
 $(3^{k+1}t_{3})d_{1}\equiv (3^{k+1}t_{3})d_{2}\equiv 0(mod\,l_{0})$, $3^{k+1}t_{3}\equiv 0(mod\,l_{0})$. Then from \eqref{T1:3}
 $C(2^{3^{k}t_{3}}-1)\equiv 0(mod\,d)$, \\$2^{3^{k}t_{3}}-1\equiv 0(mod\,d)$, $A\equiv 3\equiv 0(mod\,d)$, $d=3$. Have received the contradiction.
The equality \eqref{T1:1} is proved.\par
2) Let's prove equality \eqref{T1:2}. Let's consider the following formulas
\begin{eqnarray}
\textcolor[rgb]{0.00,0.00,0.50}{A_{1}=2^{2\cdot 3^{k_{1}}t_{1}} + 2^{3^{k_{1}}t_{1}} + 1,}
&\textcolor[rgb]{0.00,0.00,0.50}{A_{1}(2^{3^{k_{1}}t_{1}}-1)=2^{3^{k_{1}+1}t_{1}}-1,}\notag\\
\textcolor[rgb]{0.00,0.00,0.50}{A_{2}=2^{2\cdot 3^{k_{2}}t_{2}} + 2^{3^{k_{2}}t_{2}} + 1,}
&\textcolor[rgb]{0.00,0.00,0.50}{A_{2}(2^{3^{k_{2}}t_{2}}-1)=2^{3^{k_{2}+1}t_{2}}-1,}\label{T1:8}
\end{eqnarray}
where $k_{1}\neq k_{2}$, $k_{1}<k_{2}$.\par
Let's assume, that $q>1$ is a prime number such, that\\ $A_{1}\equiv A_{2}\equiv 0(mod\,q)$, $2^{3^{k_{1}+1}t_{1}}\equiv 2^{3^{k_{2}+1}t_{2}}\equiv 1(mod\,q)$.
Let $l_{0}>1$ is an index, to which the number $2$ belongs on the module $q$, $2^{l_{0}}-1\equiv 0(mod\,q)$. Then
 $3^{k_{1}+1}t_{1}\equiv 3^{k_{2}+1}t_{2}\equiv 0(mod\,l_{0})$. Since $k_{1}+1\leq k_{2}$, then\\ $3^{k_{2}}t_{2}\equiv 0(mod\,l_{0})$,
 $A_{2}\equiv 3\equiv 0(mod\,q)$, $q=3$. Have received the contradiction. The theorem 1 is proved.
\end{proof}
\begin{corollary}
If $n=3^{k}t$, where $k\geq1$ is integer, $t>1$ is an odd number, $t\not\equiv0(mod\,3)$, then
\[
\textcolor[rgb]{0.00,0.00,0.50}{L_{1}(3^{k}t)\equiv 0(mod\,L_{1}(3^{k}))}
\]
is always composite number.
\end{corollary}
The summary of the received results concerning numbers $L_{1}$.
\begin{theorem}[About numbers $L_{1}$]\label{T:T2}
For numbers $L_{1}(n) = 2^{2n} + 2^{n} + 1$ the following statements are fair :\par
 (1) The prime numbers $L_{1}(1)=7$, $L_{1}(3)=73$, $L_{1}(9)=262657$ are known.\par
 (2) If $n\ne3^{k}$, where $k\geq0$ is an integer, then $L_{1}(n)$ is a composite number.
\end{theorem}
From the identity $2^{3^{k+1}}-1=(2^{3^{k}}-1)[(2^{3^{k}})^{2}+2^{3^{k}}+1]$ the following equality is received for numbers $L_{1}$
\begin{equation}
\textcolor[rgb]{0.00,0.00,0.50}{2^{3^{k+1}}-1=L_{1}(1)\cdot L_{1}(3)\cdot\ldots\cdot L_{1}(3^{k}),}
\end{equation}
where $k\geq0$ is integer.\par
From the statement \eqref{T1:2} of the theorem 1 the property follows
\begin{equation}
\textcolor[rgb]{0.00,0.00,0.50}{L_{1}(3^{i})\bot L_{1}(3^{j})}
\end{equation}
at $i\ne j$.
\paragraph{\textbf{\textsection 3}.} The numbers $L_{1}$ are not free from squares, that is confirmed by the following examples
$L_{1}(7)\equiv 0(mod\,7^{2})$, $L_{1}(104)\equiv 0(mod\,13^{2})$, \\$L_{1}(114)\equiv 0(mod\,19^{2})$.
The following theorem takes place:
\begin{theorem}[]
Let $k\geq0$, $n\not\equiv0(mod\,3)$ are integers. Then the comparison is fair
\begin{equation}
\textcolor[rgb]{0.00,0.00,0.50}{L_{1}(7^{k}n)\equiv0(mod\,7^{k+1}).}
\end{equation}
\end{theorem}
\begin{proof}
Let's consider number $B=2^{7^{k}n}-1$. Since $n\not\equiv0(mod\,3)$, then\\ $B\not\equiv0(mod\,7)$.\par
Let $A=B\cdot L_{1}(7^{k}n)=2^{7^{k}3n}-1$. Let's prove by induction on $k\geq0$, that
\begin{equation}
\textcolor[rgb]{0.00,0.00,0.50}{2^{7^{k}3n}-1\equiv0(mod\,7^{k+1}).}\label{T3:1}
\end{equation}
The case $k=0$ is obvious. Let's make the inductive assumption, that for $k\leq m-1$ the comparison \eqref{T3:1}  is fair.
Let's consider expression
\begin{multline}
\notag \textcolor[rgb]{0.00,0.00,0.50}{2^{7^{m}3n}-1=(2^{7^{m-1}3n})^{7}-1=}\\
\textcolor[rgb]{0.00,0.00,0.50}{=(2^{7^{m-1}3n}-1)\cdot\sum_{i=0}^6(2^{7^{m-1}3n})^{i}\equiv0(mod\,7^{m}\cdot7)\equiv0(mod\,7^{m+1}).}
\end{multline}
\end{proof}

{\centering\section{Numbers $L_{2}(n) = 2^{2n} + 2^{n} - 1$}\par}

First five prime numbers $L_{2}(1)=5$, $L_{2}(2)=19$, $L_{2}(3)=71$, \\$L_{2}(4)=271$, $L_{2}(6)=4159$.
From \textbf{the statement 1} validity of the \\comparisons follows
\[
L_{2}(4k+1)\equiv0(mod\,5),\  L_{2}(10k+7)\equiv L_{2}(10k+8)\equiv0(mod\,11),
\]
where $k\geq0$ is integer.\par
The numbers $L_{2}(n)$ are not free from squares $L_{2}(68)\equiv 0(mod\,11^{3})$,\\ $L_{2}(97)\equiv 0(mod\,11^{2})$.\par
The author has checked up the following worthy to attention facts for numbers $L_{2}(n)$.\par
1) For prime numbers $p\leq5003$ the prime numbers $L_{2}(p) = 2^{2p} + 2^{p} - 1$ exist only for $p=2$,
$L_{2}(2)=19$; $p=3$, $L_{2}(3)=71$; $p=379$.\par
2) If we consider numbers $L_{2}(2^{n}) = 2^{2^{n+1}} + 2^{2^{n}} - 1$, then prime numbers $L_{2}(2^{n})$ for $n\leq17$
exist at $n=1$, $L_{2}(2)=19$; $n=2$, $L_{2}(4)=271$; $n=4$, $L_{2}(16)=4295032831$.\par

{\centering\section{Numbers $L_{3}(n) = 2^{2n} - 2^{n} + 1$}\par}

Trivial property of numbers $L_{3}(n)$ : if $n>0$ is an even number, then
\begin{equation}
\textcolor[rgb]{0.00,0.00,0.50}{L_{3}(n) \equiv 1\not\equiv 0(mod\,3),\label{L3:1}}
\end{equation}
if $n>0$ is an odd number, then
\begin{equation}
\textcolor[rgb]{0.00,0.00,0.50}{L_{3}(n) \equiv 0(mod\,3).\tag{\ensuremath{\ref{L3:1}'}}\label{L3:2}}
\end{equation}
From \textbf{the statement 1} validity of comparisons follows
\[
L_{3}(2(6k+1))\equiv L_{3}(2(6k+5))\equiv0(mod\,13),
\]
where $k\geq0$ is integer.\par
The numbers $L_{3}(n)$ are not free from squares \\$L_{3}(2\cdot13)\equiv L_{3}(10\cdot13)\equiv0(mod\,13^{2})$,
$L_{3}(3\cdot19)\equiv 0(mod\,19^{2})$.\par
In the following theorem the interesting property of $gcd$ for numbers $L_{3}$ is proved.
\begin{theorem}[$gcd$ of numbers $L_{3}$]\label{T:T4}
Let $m\geq0$, $m_{1}\geq0$, $m_{2}\geq0$, $n>0$, \\$n_{1}>0$, $n_{2}>0$ are integers;
$t_{1}\geq1$, $t_{2}\geq1$ are odd numbers, \\$t_{1}\not\equiv0(mod\,3)$, $t_{2}\not\equiv0(mod\,3)$.
Then the statements are fair :
\begin{gather}
\label{T4:1}\textcolor[rgb]{0.00,0.00,0.50}{gcd(L_{3}(3^{m}2^{n}t_{1}),L_{3}(3^{m}2^{n}t_{2}))=
L_{3}(3^{m}2^{n}gcd(t_{1},t_{2})).}\\
\label{T4:2}\textcolor[rgb]{0.00,0.00,0.50}{gcd(L_{3}(3^{m_{1}}2^{n_{1}}t_{1}),L_{3}(3^{m_{2}}2^{n_{2}}t_{2})) = 1}
\end{gather}
for $m_{1}\neq m_{2}$ or $n_{1}\neq n_{2}$.
\end{theorem}
\begin{proof}
Let $t_{3}=gcd(t_{1},t_{2})$, $t_{1}=t_{3}d_{1}$, $t_{2}=t_{3}d_{2}$,
where \\$(d_{1},d_{2})=1$, $d_{1}\geq1$, $d_{2}\geq1$ are odd numbers.\par
1) Let's prove equality \eqref{T4:1}. Let's consider the following formulas
\begin{gather}
\textcolor[rgb]{0.00,0.00,0.50}{A=L_{3}(3^{m}2^{n}t_{1})=2^{2\cdot 3^{m}2^{n}t_{1}} - 2^{3^{m}2^{n}t_{1}} + 1,}\notag\\
\textcolor[rgb]{0.00,0.00,0.50}{A(2^{3^{m}2^{n}t_{1}}+1)=2^{3^{m+1}2^{n}t_{3}d_{1}}+1,}\notag\\
\textcolor[rgb]{0.00,0.00,0.50}{B=L_{3}(3^{m}2^{n}t_{2})=2^{2\cdot 3^{m}2^{n}t_{2}} - 2^{3^{m}2^{n}t_{2}} + 1,}\label{T4:3}\\
\textcolor[rgb]{0.00,0.00,0.50}{B(2^{3^{m}2^{n}t_{2}}+1)=2^{3^{m+1}2^{n}t_{3}d_{2}}+1,}\notag\\
\textcolor[rgb]{0.00,0.00,0.50}{C=L_{3}(3^{m}2^{n}t_{3})=2^{2\cdot 3^{m}2^{n}t_{3}} - 2^{3^{m}2^{n}t_{3}} + 1},\notag\\
\textcolor[rgb]{0.00,0.00,0.50}{C(2^{3^{m}2^{n}t_{3}}+1)=2^{3^{m+1}2^{n}t_{3}}+1.}\notag
\end{gather}
Then the formulas are fair
\begin{multline}
\textcolor[rgb]{0.00,0.00,0.50}{A(2^{3^{m}2^{n}t_{1}}+1)=C(2^{3^{m}2^{n}t_{3}}+1)[2^{3^{m+1}2^{n}t_{3}(d_{1}-1)} -}\\
\textcolor[rgb]{0.00,0.00,0.50}{-2^{3^{m+1}2^{n}t_{3}(d_{1}-2)} + \ldots}\\
\textcolor[rgb]{0.00,0.00,0.50}{\ldots + 2^{3^{m+1}2^{n}t_{3}2} - 2^{3^{m+1}2^{n}t_{3}} + 1],}\label{T4:4}
\end{multline}
\begin{multline}
\textcolor[rgb]{0.00,0.00,0.50}{B(2^{3^{m}2^{n}t_{2}}+1)=C(2^{3^{m}2^{n}t_{3}}+1)[2^{3^{m+1}2^{n}t_{3}(d_{2}-1)} -}\\
\textcolor[rgb]{0.00,0.00,0.50}{-2^{3^{m+1}2^{n}t_{3}(d_{2}-2)} + \ldots}\\
\textcolor[rgb]{0.00,0.00,0.50}{\ldots + 2^{3^{m+1}2^{n}t_{3}2} - 2^{3^{m+1}2^{n}t_{3}} + 1].}\label{T4:5}
\end{multline}
Let $q>1$ is prime number such, that
\begin{equation}
\textcolor[rgb]{0.00,0.00,0.50}{2^{3^{m}2^{n}t_{1}}+1\equiv 0(mod\,q),\ C\equiv 0(mod\,q).}\label{T4:6}
\end{equation}
Then from \eqref{T4:3} and \eqref{T4:6} the comparisons follow
\[
\textcolor[rgb]{0.00,0.00,0.50}{2^{3^{m}2^{n}t_{3}d_{1}}+1\equiv 0(mod\,q),\ 2^{3^{m+1}2^{n}t_{3}}+1\equiv 0(mod\,q).}
\]
Let's consider number $b=2^{3^{m}2^{n}t_{3}}$. If $b\equiv 1(mod\,q)$, then \\$C=b^{2}-b+1\equiv 1\not\equiv0(mod\,q)$.
If $b\equiv -1(mod\,q)$, then \\$C\equiv 3\equiv 0(mod\,q)$, $q=3$, but as $n>0$, that from \eqref{L3:1} follows, that\\
$C\not\equiv 0(mod\,3)$. Have received the contradiction. Thus, $b^{2}\not\equiv 1(mod\,q)$.\par
Let $l_{0}>1$ is an index, to which the number $b^{2}$ belongs on the module $q$, $(b^{2})^{l_{0}}\equiv 1(mod\,q)$.
As $(b^{2})^{d_{1}}\equiv (b^{2})^{3}\equiv 1(mod\,q)$, that\\
$d_{1}\equiv 3\equiv 0(mod\,l_{0})$, $d_{1}\equiv 0(mod\,3)$. Have received the contradiction.\\
Have proved, that $gcd(A,B)\equiv 0(mod\,C)$.\par Let's assume, that there is an integer $d>1$ such, that
\begin{equation}
\textcolor[rgb]{0.00,0.00,0.50}{A\equiv B\equiv 0(mod\,d),\ but\ gcd(C,d)=1.}\label{T4:7}
\end{equation}
Then it follows from \eqref{T4:3} and \eqref{T4:7}, that
\[
\textcolor[rgb]{0.00,0.00,0.50}{2^{3^{m+1}2^{n+1}t_{3}d_{1}}-1\equiv 0(mod\,d),\
2^{3^{m+1}2^{n+1}t_{3}d_{2}}-1\equiv 0(mod\,d).}
\]
Then $2^{3^{m+1}2^{n+1}t_{3}}-1\equiv 0(mod\,d)$, i.e.
\begin{equation}
\textcolor[rgb]{0.00,0.00,0.50}{(2^{3^{m+1}2^{n}t_{3}}+1)\cdot(2^{3^{m+1}2^{n}t_{3}}-1)\equiv 0(mod\,d).}\label{T4:8}
\end{equation}
If $gcd(2^{3^{m+1}2^{n}t_{3}}+1,d)=d_{0}>1$, then it follows from \eqref{T4:3}, that \\$2^{3^{m}2^{n}t_{3}}\equiv -1(mod\,d_{0})$,
$A\equiv 3\equiv 0(mod\,d_{0})$, $d_{0}=3$. Have received the contradiction. Then from \eqref{T4:8}, \eqref{T4:7} and
\eqref{T4:3} the comparisons follow\\ $(2^{3^{m+1}2^{n}t_{3}}-1)\equiv 0(mod\,d)$, $2^{3^{m+1}2^{n}t_{3}d_{1}}+1\equiv 0\equiv 2(mod\,d)$.
Have received the contradiction, since $d>1$ is odd number. The equality \eqref{T4:1} is proved.\par
2) Let's prove equality \eqref{T4:2}. Let's consider the following formulas
\begin{gather}
\textcolor[rgb]{0.00,0.00,0.50}{A_{1}=2^{2\cdot 3^{m_{1}}2^{n_{1}}t_{1}} - 2^{3^{m_{1}}2^{n_{1}}t_{1}} + 1,}\notag\\
\textcolor[rgb]{0.00,0.00,0.50}{A_{1}(2^{3^{m_{1}}2^{n_{1}}t_{1}}+1)=2^{3^{m_{1}+1}2^{n_{1}}t_{3}d_{1}}+1,}\label{T4:9}\\
\textcolor[rgb]{0.00,0.00,0.50}{A_{2}=2^{2\cdot 3^{m_{2}}2^{n_{2}}t_{2}} - 2^{3^{m_{2}}2^{n_{2}}t_{2}} + 1,}\notag\\
\textcolor[rgb]{0.00,0.00,0.50}{A_{2}(2^{3^{m_{2}}2^{n_{2}}t_{2}}+1)=2^{3^{m_{2}+1}2^{n_{2}}t_{3}d_{2}}+1.}\notag
\end{gather}
Let's assume, that $q>1$ is prime number such, that\\
$A_{1}\equiv A_{2}\equiv 0(mod\,q)$. Then $2^{3^{m_{1}+1}2^{n_{1}}t_{3}d_{1}}+1\equiv 2^{3^{m_{2}+1}2^{n_{2}}t_{3}d_{2}}+1\equiv 0(mod\,q)$.
Let $l_{0}>1$ is an index, to which the number $2$ belongs on the module $q$. Then
 $3^{m_{1}+1}2^{n_{1}+1}t_{3}d_{1}\equiv 3^{m_{2}+1}2^{n_{2}+1}t_{3}d_{2}\equiv 0(mod\,l_{0})$, i.e.
\begin{equation}
\textcolor[rgb]{0.00,0.00,0.50}{3^{m_{1}+1}2^{n_{1}+1}t_{3}\equiv 3^{m_{2}+1}2^{n_{2}+1}t_{3}\equiv 0(mod\,l_{0})}\label{T4:10}.
\end{equation}\par
$2^{a}$) Let's assume, that $m_{1}<m_{2}$. Then from \eqref{T4:10} we receive comparisons
\[
3^{m_{2}}2^{n_{2}+1}t_{3}\equiv 0(mod\,l_{0}),\ (2^{3^{m_{2}}2^{n_{2}}t_{3}}+1)\cdot(2^{3^{m_{2}}2^{n_{2}}t_{3}}-1)\equiv 0(mod\,q).
\]
From the last comparison either $A_{2}\equiv 3\equiv 0(mod\,q)$, $q=3$,
or \\$A_{2}\equiv 1\not\equiv0(mod\,q)$ follows. Have received the contradiction.\par
$2^{b}$) Let's assume, that $m_{1}=m_{2}$, $n_{1}<n_{2}$. Then from \eqref{T4:10}, \eqref{T4:9} we receive comparisons
\[
3^{m_{2}+1}2^{n_{2}}t_{3}\equiv 0(mod\,l_{0}),\ (2^{3^{m_{2}+1}2^{n_{2}}t_{3}d_{2}}+1)\equiv 0\equiv 2(mod\,q).
\]
Have received the contradiction. The theorem 4 is proved.
\end{proof}

\begin{corollary}
Let $n\geq1$ is an integer, $t>1$ is an odd number. Then the statements are fair :\par
(1) If $t\not\equiv0(mod\,3)$, then
\[
\textcolor[rgb]{0.00,0.00,0.50}{L_{3}(2^{n}t)\equiv 0(mod\,L_{3}(2^{n})).}
\]
Besides, $1<L_{3}(2^{n})<L_{3}(2^{n}t)$, where $L_{3}(2^{n}t)$ is composite number.\par
(2) If $t\equiv0(mod\,3)$, then
\[
\textcolor[rgb]{0.00,0.00,0.50}{gcd(L_{3}(2^{n}t),L_{3}(2^{n}))=1.}
\]
\end{corollary}

\begin{corollary}
Let $m\geq0$, $n>0$ are integers; $t>1$ is an odd number, $t\not\equiv0(mod\,3)$.
Then number $L_{3}(3^{m}2^{n}t)$ is always composite number.
\end{corollary}\par
The summary of the received results concerning composite numbers $L_{3}$.\par
\begin{theorem}[About numbers $L_{3}$]\label{T:T5}
For numbers $L_{3}$ the statements are fair :\par
(1) If $n\ne 3^{m}2^{n}$, where $m\geq0$, $n\geq0$ are integers, then number $L_{3}(n)$ is composite number.\par
(2) Prime numbers $L_{3}(1)=L_{3}(2^{0})=3$, $L_{3}(2)=L_{3}(2^{1})=13$,
$L_{3}(4)=L_{3}(2^{2})=241$, $L_{3}(32)=L_{3}(2^{5})=18446744069414584321$ are known.
\end{theorem}
\begin{proof}
\end{proof}\par
For numbers $L_{3}(3^{m}2^{n})$, where $m\geq0$, $n>0$ are integers, the author has carried out the following check :\par
1) Numbers $L_{3}(2^{n})$ at $6\leq n \leq 15$ is composite;\par
2) Numbers $L_{3}(2\cdot3^{m})$ at $1\leq m \leq 8$ is composite;\par
3) Numbers $L_{3}(3\cdot2^{n})$ at $1\leq n \leq 12$ is composite;\par
4) Numbers $L_{3}(3^{2}\cdot2^{n})$ at $1\leq n \leq 11$ is composite;\par
5) Numbers $L_{3}(3^{3}\cdot2^{n})$ at $1\leq n \leq 9$ is composite.\par

{\centering\section{Numbers $L_{4}(n) = 2^{2n} - 2^{n} - 1$}\par}

From \textbf{the statement 1} validity of comparisons follows
\[
L_{4}(4k+3)\equiv0(mod\,5),\  L_{4}(10k+2)\equiv L_{4}(10k+3)\equiv0(mod\,11),
\]
where $k\geq0$ is any integer.\par
The numbers $L_{4}(n)$ are not free from squares, since
\begin{gather}
 L_{4}(13)\equiv L_{4}(42)\equiv L_{4}(123)\equiv0(mod\,11^{2})\notag,\\
 L_{4}(52)\equiv L_{4}(119)\equiv 0(mod\,19^{2})\notag.
\end{gather}\par
Prime numbers $L_{4}(n)$ and $L_{4}(n+1)$ are named \textbf{the prime $L_{4}$ number-twins}.
\textbf{The author has found 4 pairs of the prime $L_{4}$ number-twins up to $n\leq603$}, namely
\begin{center}
$L_{4}(1) = 1, L_{4}(2) = 11$ ; $L_{4}(4) = 239, L_{4}(5) = 991$ ;\\
$L_{4}(9) = 261631, L_{4}(10) = 1047551$ ; $L_{4}(224) , L_{4}(225)$.
\end{center}

{\centering\section{Wonderful properties of $gcd$ insularity}\par}
\begin{definition}[Insularity to $gcd$.]
Let $R_{n}\geq0$ is a sequence of integers, where $n>0$ is integer. $M$ is a subset of natural numbers.\\
Let's tell, that the sequence $R_{n}$ on set $M$ is isolated about $gcd$, if the condition is fair :
\begin{equation}
\textcolor[rgb]{0.00,0.00,0.50}{gcd(R_{n},R_{m})=R_{gcd(n,m)}}
\end{equation}
for $\forall n,m \in M$.
\end{definition}
\begin{corollary}
Let $k\geq0$ is integer, then numbers $L_{1}$ on set\par
$M_{k}=\{ 3^{k}\cdot t : t\geq1$ is an odd number, $t\not\equiv 0(mod\,3)\}$\\
are isolated about $gcd$, i.e.
\begin{equation}
\textcolor[rgb]{0.00,0.00,0.50}{gcd(L_{1}(n),L_{1}(m))=L_{1}(gcd(n,m))}\notag
\end{equation}
for $\forall n,m \in M_{k}$.
\end{corollary}

\begin{corollary}
Let $k\geq0, l>0$ are integers, then numbers $L_{3}$ on set\par
$M_{k,l}=\{ 3^{k}2^{l}\cdot t : t\geq1$ is an odd number, $t\not\equiv 0(mod\,3)\}$\\
are isolated about $gcd$, i.e.
\begin{equation}
\textcolor[rgb]{0.00,0.00,0.50}{gcd(L_{3}(n),L_{3}(m))=L_{3}(gcd(n,m))}\notag
\end{equation}
for $\forall n,m \in M_{k,l}$.
\end{corollary}


Let's consider the generalized numbers repunit - integers of the following kind
\cite{WeissteinM,WeissteinR,Dubner-1,Dubner-2,Granlund} :
\begin{equation}\label{GenRep1}
\textcolor[rgb]{0.00,0.00,0.50}{M_{n}^{(b)} = (b^{n}-1)/(b-1),}
\end{equation}
where $n\geq1$, $b\geq2$ are integers.
\begin{equation}\label{GenRep2}
\textcolor[rgb]{0.00,0.00,0.50}{M_{n}^{+(b)} = (b^{n}+1)/(b+1),}
\end{equation}
where $n\geq1$ is an odd number, $b\geq2$ is integer.\par
For the generalized numbers repunit \eqref{GenRep1} and \eqref{GenRep2} the theorem takes place.
\begin{theorem}[]\label{T:T6} Following formulas are fair :
\begin{equation}\label{e:Form1}
\textcolor[rgb]{0.00,0.00,0.50}{gcd(M_{n}^{(b)},M_{m}^{(b)}) = M_{gcd(n,m)}^{(b)}},
\end{equation}
where $n\geq1$, $b\geq2$ are integers.
\begin{equation}\label{e:Form2}
\textcolor[rgb]{0.00,0.00,0.50}{gcd(M_{n}^{+(b)},M_{m}^{+(b)}) = M_{gcd(n,m)}^{+(b)}},
\end{equation}
where $n\geq1$ is an odd number, $b\geq2$ is integer.\par
\end{theorem}
\begin{proof}\par
1) Let $(n,m)=d\geq1$, where $n=n_{1}d$, $m=m_{1}d$, $(n_{1},m_{1})=1$.
From definition \eqref{GenRep1} equalities follow
\begin{displaymath}M_{n}^{(b)} = ((b^{d})^{n_{1}}-1)/(b-1)=M_{d}^{(b)}\cdot\{b^{d(n_1-1)}+\,\ldots\,+b^{2d}+b^d+1\},\end{displaymath}
\begin{displaymath}M_{m}^{(b)} = ((b^{d})^{m_{1}}-1)/(b-1)=M_{d}^{(b)}\cdot\{b^{d(m_1-1)}+\,\ldots\,+b^{2d}+b^d+1\}.\end{displaymath}\par
Let
\begin{displaymath}A=b^{d(n_1-1)}+\,\ldots\,+b^{2d}+b^d+1,\ \ B=b^{d(m_1-1)}+\,\ldots\,+b^{2d}+b^d+1.\end{displaymath}\par
Let's assume, that $A\equiv~B\equiv0(mod\,q)$, where $q>1$ is prime number.\\
Let $b_{0}=b^{d}$. If $b_{0}\equiv1(mod\,q)$, then $n_{1}\equiv~m_{1}\equiv0(mod\,q)$. Have received the contradiction.
Hence $b_{0}\not\equiv1(mod\,q)$, then there exists an index $d_{0}>1$, to which the number $b_{0}$ belongs on
the module $q$
\begin{displaymath}(b^d)^{d_{0}}\equiv1(mod\,q).\end{displaymath}
Then $n_{1}\equiv m_{1}\equiv0(mod\,d_{0})$. Have received the contradiction.\par
2) Let $(n,m)=d\geq1$, where $n=n_{1}d$, $m=m_{1}d$ are odd numbers, $(n_{1},m_{1})=1$.
From definition \eqref{GenRep2} equalities follow
\begin{displaymath}M_{n}^{+(b)} = ((b^{d})^{n_{1}}+1)/(b+1) = \end{displaymath}
\begin{displaymath} = M_{d}^{+(b)}\cdot\{b^{d(n_1-1)}-b^{d(n_1-2)}+\,\ldots\,+b^{2d}-b^d+1\},\end{displaymath}
\begin{displaymath}M_{m}^{+(b)} = ((b^{d})^{m_{1}}+1)/(b+1) = \end{displaymath}
\begin{displaymath} = M_{d}^{+(b)}\cdot\{b^{d(m_1-1)}-b^{d(m_1-2)}+\,\ldots\,+b^{2d}-b^d+1\}.\end{displaymath}\par
Let
\begin{displaymath}A=b^{d(n_1-1)}-b^{d(n_1-2)}+\,\ldots\,+b^{2d}-b^d+1,\end{displaymath}
\begin{displaymath}B=b^{d(m_1-1)}-b^{d(m_1-2)}+\,\ldots\,+b^{2d}-b^d+1.\end{displaymath}\par
Let's assume, that $A\equiv~B\equiv0(mod\,q)$, where $q>1$ is prime number.\\
Let $b_{0}=b^{2d}$. If $b_{0}\equiv1(mod\,q)$, then either $b^{d}\equiv1(mod\,q)$, or\\ $b^{d}\equiv-1(mod\,q)$.
Then either $A\equiv1\not\equiv0(mod\,q)$, or\\ $n_{1}\equiv~m_{1}\equiv0(mod\,q)$. Have received the contradiction.

Hence $b_{0}\not\equiv1(mod\,q)$, then there exists an index $d_{0}>1$, to which the number $b_{0}$ belongs on
the module $q$
\begin{displaymath}(b^{2d})^{d_{0}}\equiv1(mod\,q).\end{displaymath}
Since $(b^{2d})^{n_{1}}\equiv (b^{2d})^{m_{1}}\equiv1(mod\,q)$, then $n_{1}\equiv m_{1}\equiv0(mod\,d_{0})$. Have received the contradiction.
\end{proof}


\begin{corollary}
Let $\mathbb{P}$ is a set of all positive integers, $\mathbb{O}$ is a set of all odd numbers.\par
(1) Numbers $M_{n}^{(b)}$ on set $\mathbb{P}$ are isolated about $gcd$, i.e.
\[
\textcolor[rgb]{0.00,0.00,0.50}{gcd(M_{n}^{(b)},M_{m}^{(b)}) = M_{gcd(n,m)}^{(b)}}
\]
for $\forall n,m \in \mathbb{P}$.\par

(2) Numbers $M_{n}^{+(b)}$ on set $\mathbb{O}$ are isolated about $gcd$, i.e.
\[
\textcolor[rgb]{0.00,0.00,0.50}{gcd(M_{n}^{+(b)},M_{m}^{+(b)}) = M_{gcd(n,m)}^{+(b)}}
\]
for $\forall n,m \in \mathbb{O}$.
\end{corollary}

{\centering\section{The open problems of numbers $L(n)$}\par}

Author offers some open problems, as the unsolved tasks concerning numbers $L(n) = 2^{2n}\pm2^n\pm1$.
\begin{problem}[]
\textcolor[rgb]{0.00,0.00,0.50}{Whether there are prime numbers $L_{1}(3^{k})$ for $k>2$} ?
\end{problem}
The author has checked up, that the numbers $L_{1}(3^{k})$ for $k=3,4,5,6,7,\\8,9,10$ are composite !
\begin{problem}[]
\textcolor[rgb]{0.00,0.00,0.50}{Whether there are infinitely many prime numbers $L_{2}(p)$, where $p$ is prime number} ?
\end{problem}
\begin{problem}[]
\textcolor[rgb]{0.00,0.00,0.50}{Whether there are prime numbers $L_{2}(2^{n})$ for $n>17$} ?
\end{problem}
\begin{problem}[]
\textcolor[rgb]{0.00,0.00,0.50}{Whether there are prime numbers $L_{3}(2^{n})$ for $n>5$} ?
\end{problem}
\begin{problem}[]
\textcolor[rgb]{0.00,0.00,0.50}{Whether there are prime numbers $L_{3}(3^{m}\cdot2^{n})$ for $m\geq1$, \\$n>1$} ?
\end{problem}
\begin{problem}[]
\textcolor[rgb]{0.00,0.00,0.50}{Whether there are infinitely many prime numbers-twins $L_{4}(n)$, \\$L_{4}(n+1)$, where $n\geq1$} ?
\end{problem}

{\centering\section{Conclusion "The concrete theory of numbers"}\par}

It is necessary to explain the title of article "The concrete theory of numbers".
Having had a look in the Wladimir Igorewitsch Arnold foreword "From Fibonacci up to Erd\capitaldieresis{o}s" to the remarkable
book of Ronald Graham, Donald Knuth and Oren Patashnik "The Concrete mathematics"\cite{Grekhem},
it is possible to answer the question :
\textbf{what is "the concrete theory of numbers" ?}\par
\textbf{The theories come and leave, but natural series of numbers remains} and constantly generates new
complicated problems.
The new theories are again created for their decision. The process cannot be stopped \textbf{:)}\par
The concrete theory of numbers created by titanic efforts of Pierre de Fermat and Leonhard Euler is an art to solve riddles of a natural series of numbers.
It is enough to look to the tasks list \cite{Zadachi Ferma}, which has been put and decided by Pierre de Fermat,
or to get acquainted with tasks, which has been decided, investigated
and propagandized by Waclaw Sierpinski \cite{Serpinskii},to be convinced -
\textbf{a natural series of numbers doesn't drowse, it is always ready to a human challenge!}

{\centering\section{Acknowledgement of gratitude}\par}

The author expresses the deep gratitude to the creators of the calculator for number-theoretic researches. --
GR/PARI CALCULATOR is free software \cite{GP/PARI}.
GR/PARI CALCULATOR has opened to the author a door in the world of wonderful number-theoretic opportunities !\par


\par

---------------------------------------------------------------------\par
Institute of Thermophysics, Siberian Branch of RAS \par
Lavrentyev Ave., 1, Novosibirsk, 630090, Russia \par
E-mail: tarasov@itp.nsc.ru \par
---------------------------------------------------------------------\par
Boris Vladimirovich Tarasov,
independent researcher.\par
Primary E-mail Address: tarasov-b@mail.ru
\end{document}